\documentclass[a4paper,reqno,10pt,oneside]{amsart}

\usepackage{amsmath,amsthm,amsfonts,amssymb}
\usepackage[latin1]{inputenc} % letras acentuadas
\usepackage{todonotes}
\usepackage{tikz}
\usetikzlibrary{shapes,arrows}
\usepackage[numbers]{natbib}
\usepackage{a4wide}
\usepackage{graphicx}

\newtheorem{theorem}{Theorem}[section]                                          
\newtheorem{prop}[theorem]{Proposition}

\newtheorem{definition}[theorem]{Definition}
\newtheorem{remark}{Remark}[section]

\makeatletter
\@addtoreset{equation}{section}
\makeatother

\renewcommand\thetable{\thesection.\@arabic\c@table}

\bibdata{prob}
\bibstyle{alpha} 

\title[A two-stage innovation diffusion model]{A stochastic two-stage innovation diffusion model\\ on a lattice}

\author{Cristian F. Coletti, Karina B. E. de Oliveira and Pablo M. Rodriguez}

\date{}

\address{
\newline
UFABC - Centro de Matem\'atica, Computa\c{c}\~ao e Cogni\c{c}\~ao
\newline
Avenida dos Estados, 5001- Bangu - Santo André - São Paulo, Brasil
\newline
e-mail:  cristian.coletti@ufabc.edu.br
\newline
\newline
USP, Instituto de Ciências Matemáticas e de Computação
\newline  
Av. Trabalhador são-carlense 400 - Centro, CEP 13560-970, São Carlos, SP, Brasil
\newline
e-mail: emboaba@icmc.usp.br
\newline
e-mail:  pablor@icmc.usp.br
}

\begin{document}
  
%\maketitle

\begin{abstract}
We propose a stochastic model describing a process of awareness, evaluation and decision-making by agents on the $d$-dimensional integer lattice. Each agent may be in any of the three states belonging to the set $\{0, 1, 2\}$. In this model $0$ stands for ignorants, $1$ for aware and $2$ for adopters.  Aware and adopters inform its nearest ignorant neighbors about a new product innovation at rate $\lambda$. At rate $\alpha$ an agent in aware state becomes an adopter due to the influence of adopters neighbors. Finally, aware and adopters forget the information about the new product, {\color{black}thus becoming ignorant}, at rate one. Our purpose is to analyze the influence of the parameters on the qualitative behavior of the process. We obtain sufficient conditions under which the innovation diffusion (and adoption) either becomes extinct or propagates through the population with positive probability.
\end{abstract}

\maketitle

%%%%%%%%%%%%%%%%%%%%%%%%%%%%%%%%%%%%%%%%%%%%%%%%%%%%%%%%%%%%%%% S: INTRODUÇÃO
%%%%%%%%%%%%%%%%%%%%%%%%%%%%%%%%%%%%%%%%%%%

\section{Introduction} 

\label{sec:intro}

The study of dissemination of information (news, innovations, rumors) nowadays is extremely important and become the focus of much research --- see for instance \cite{agliari/2010,agliari/2006,arruda/2015,journals/siamads/McCullenRBFG13,raey,CPC:9370507,comets/gallesco/popov/vachkovskaia/2013,isham/harden/nekovee/2010,lebensztayn/machado/rodriguez/2011a,lebensztayn/machado/rodriguez/2011b}, and references therein. The purpose of this work is to analyze, by means of a mathematical model, the diffusion of an innovation which has been defined as ``the process by which an innovation is communicated through certain channels over time among the members of a social system,'' Rogers \cite{rogers1983diffusion}. According to Rogers \cite{rogers1983diffusion} the adoption of an innovation occurs through a process which could be divided into 5 stages: {\color{black}
\begin{enumerate}	
\item	 knowledge (the individual is introduced to an innovation), 
\item persuasion (the individual forms a favorable or unfavorable attitude toward to adopt or reject the innovation), 
\item decision (the individual engages in activities that lead to a choice to adopt the innovation or reject it), 
\item implementation (the individual puts an innovation into use), and 
\item confirmation (the individual seeks reinforcement on a decision about an innovation but may reverse this decision if exposed to conflicting messages about the innovation). 
\end{enumerate}}
	One of the first mathematical models describing the diffusion of a new product was  introduced by Bass \cite{Bass1969} , who inspired by Rogers' work \cite{rogers1983diffusion},  assumes that sales of a new product are essentially stimulated by word-of-mouth from satisfied customers.
	
	The Bass model assumes that adopters of an innovation are subdivided into two groups: innovators and imitators. The decision of {\color{black}the innovators} is influenced only by mass media or other external influences. Such individuals decide to adopt an innovation independently of the decisions of other individuals in a social system. On the other hand, the decision of {\color{black}the imitators} is influenced only by word-of-mouth communication or other influence from those who {\color{black}have} already used the product (internal influences). More precisely, if $n$ is the total number of agents who will eventually use the product and {\color{black}$A(t)$} is the number of adopters at time $t$, the Bass model is given by the differential equation
	
$$
\frac{d A(t)}{d t} = \left(p + \frac{q}{n} A(t)\right)(n - A(t)),
$$

where parameter $p$ is the coefficient of innovation and parameter $q$ is the coefficient of imitation. Thenceforth, the Bass model has become one of the most important and widely used models in marketing research. 

There are many variants to examine in order to extend this model to more realistic scenarios. We refer the reader to \cite{Bass2004} for a discussion of some extensions and examples of applications of such model. In this paper we propose a stochastic model which takes into account essentially two kind of modifications. The first one goes in the direction of generalizing the dynamic of the process. In other words, we consider an innovation diffusion model with stage structure as suggested by Rogers \cite{rogers1983diffusion}. In this context, Wang et al. \cite{Wang2006129} extend the Bass model by introducing a two-stage structure, namely, the stage of awareness of information and the stage of decision-making. However, this model maintains the assumption that the population is homogeneously mixed. Indeed, they consider a system of differential equations and provide sufficient conditions for the success of the innovation's diffusion. Also, the authors introduce a model with a time delay {\color{black}for which} they prove the existence of stability switches.
 
 An alternative generalization is to consider a model which still exhibits a simple dynamic, as in the Bass model, but which is defined for a population with a certain neighborhood structure. In this case, McCullen et al. \cite{journals/siamads/McCullenRBFG13} propose a variant of the Bass model on complex networks where the decision to adopt the innovation takes into account not only individual preferences but also whether or not an individual's social circle has adopted it. The authors examined, using simulations, the number of neighbours needed to induce uptake and the probability of induced uptake in random networks. 

The contribution of this paper is to define and study a stochastic process that incorporates both types of generalizations in a unique model. More precisely, we consider a spatial stochastic model for a two-stage innovation diffusion. Our model may be seen as a continuous time extension of the cellular automata model examined through simulations by \cite{Goldenberg2004}, and thus, under some conditions, it may be seen as a spatial version of the Bass model. Therefore, our model may contribute to improve our understanding of aggregate behavior in innovation diffusion. 

The paper is organized as follows. In section \ref{sec:model} we introduce the model, state the main results of this work and provide a graphical construction of the stochastic process considered. In section \ref{sec:proofs} we prove the extinction of innovation awareness. We also prove the survival and extinction of the innovation adoption. Section \ref{concluding} is devoted to concluding remarks.

%%%%%%%%%%%%%%%%%%%%%%%%%%%%%%%%%%%%%%%%%%%%%%%%%%%%
%%%%%%%% S: The Model: Graphical Construction and Results
%%%%%%%%%%%%%%%%%%%%%%%%%%%%%%%%%%%%%%%%%%%%%%%%%%%%

\section{The Model: Graphical Construction and Results} 
\label{sec:model}

We consider a continuous-time Markov process $(\eta_t)_{t\geq 0}$ with {\color{black}state} space $\mathcal{S}=\{0,1,2\}^{\mathbb{Z}^d}$, i.e. at time $t$ the state of the process is some function $\eta_t: \mathbb{Z}^d \longrightarrow \{0,1,2\}$. We assume that each site $x \in \mathbb{Z}^d$ represents an agent, which is said to be an ignorant if $\eta(x)=0,$ an aware if $\eta(x)=1$ and an adopter if $\eta(x)=2.$ Ignorants are those who do not know about the innovation, aware are those who know about the innovation but they have not adopted it yet. Finally, adopters are those who have already adopted the innovation. Then, if the system is in configuration $\eta \in \mathcal{S},$ the state of site $x$ changes according to the following transition rates

\begin{equation}\label{rates}
\begin{array}{rclc}
0 & \rightarrow & 1, & \hspace{.5cm} \lambda \, (n_1(x,\eta) + n_2(x,\eta)),\\
1 & \rightarrow & 2, & \hspace{.5cm}\alpha \, n_2(x,\eta),\\
1 & \rightarrow & 0, & \hspace{.5cm}1,\\ 
2 & \rightarrow & 0, & \hspace{.5cm}1,\\ 
\end{array}
\end{equation}

\noindent where $$n_i(x,\eta)= \sum_{||x-y||=1} 1\{\eta(y)=i\}$$ 
is the number of nearest neighbors of site $x$ in state $i$ for the configuration $\eta$, for $i=1,2.$ Formally, {\color{black}(1)} means that if the site $x$ is in state, say, $0$ at time {\color{black}$t$} then the probability that it will be in state $1$ at time $t+h$, for $h$ small, is $(\lambda n_{1}(x,\eta) + \lambda n_{2}(x,\eta) )h + o(h)$, where $o(h)$ is such that $\lim_{h\to 0} o(h)/h = 0$. We call the Markov process $(\eta_t)_{t\geq 0}$ thus obtained the {\it{innovation process}} on $\mathbb{Z}^d$ with rates $\lambda$ and $\alpha$.  {\color{black} When there are no agents in state $2$, we recover the well known $d$-dimensional contact process with parameters $\lambda$ and $1$. We refer the reader to \cite[Part I]{opac-b1095541} for more details.}

In the context of an innovation diffusion scenery, {\color{black}(1)} represents the different transitions assumed in the model. The first one is related to the contact between an ignorant and an aware or adopter agent, which implies that the ignorant becomes an aware agent. The second transition represents the result of an interaction between an aware and an adopter. In this case, we assume that the adopter persuades the aware about acquiring the innovation and therefore the aware becomes an adopter. Finally, we assume that agents informed about the innovation forget about it at rate $1$. We point out that the way our model is defined, it may be used as a basis for the construction of theoretical models  that generalize the model introduced by Bass \cite{Bass1969}. Indeed, our model may be seen as a continuous time version of the model introduced by \cite{Goldenberg2004}, which assumes a stage structure such as the one proposed by Wang \cite{Wang2006129}.

%%%%%%%%%%%%%%%%%%%%%%%%%%%%%%%%%%%%%%%%%%%%%%%%%%%%
%%%%%%%% SS: Harris' Graphical Construction
%%%%%%%%%%%%%%%%%%%%%%%%%%%%%%%%%%%%%%%%%%%%%%%%%%%%

\subsection{Harris' Graphical Construction}

The well known Harris' graphical construction \cite{harris} is a powerful tool to deal with interacting particle systems. The main idea behind it, is the construction of a version of the spatial stochastic process by mean of collections of independent Poisson processes. This technique allows us to obtain interesting results by comparison (coupling) with other processes like the contact process and oriented percolation models. 

In order to obtain the graphical construction for our model, we consider a collection of independent Poisson processes denoted by $\{N^{x,y}_\lambda , N^{x,y}_\alpha , D^x : x,y \in \mathbb{Z}^d , \|x-y \|=1\}$. We assume that the Poisson processes $N^{x,y}_\lambda , N^{x,y}_\alpha$ and $D^x$ have intensities $\lambda, \alpha$ and $1$, respectively. At each arrival time of the process $N^{x,y}_\lambda$ if $x$ and $y$ are in states $1$ or $2$, and $0$, respectively then the state of $y$ is updated to state $1$. On the other hand, at each arrival time of the process $N^{x,y}_\alpha$ if sites $x$ and $y$ are in states $2$ and $1$ respectively then, the state of $y$ changes to state $2$. A last transition is obtained at an arrival time of the process $D_x$, in which case the state of site $x$ becomes $0$, provided it was in state $1$ or $2$. In this way we obtain a version of the spatial stochastic innovation process with the rates given by {\color{black}(1)}. Fig. 1 shows a possible realization of the model by means of the graphical construction. We observe that to construct the process inside a finite space-time box it is sufficient to consider the Poisson arrival times inside that box. For further details on the graphical construction {\color{black}we refer the reader to Durrett \cite[Section 2]{Durret}.}

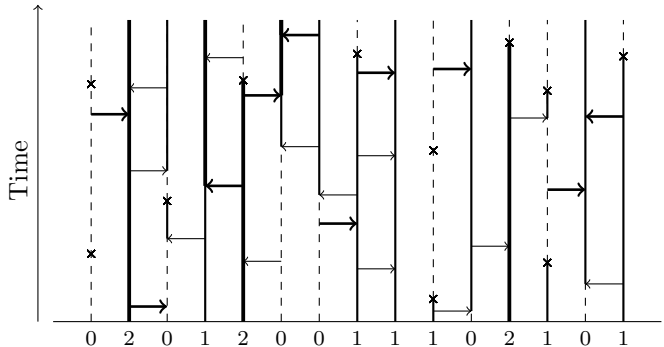
\begin{figure}[h]\label{harris}
\begin{center}
\begin{tikzpicture}
\draw (-4,-2) -- (4,-2);

% linha de tempo
\draw[->] (-4.2,-2) -- (-4.2,2.2);
\draw (-4.7,0) node[below, rotate=90] {Time};

% marcas informantes/spreaders
\draw[->] (-3,0) -- (-2.5,0);
\draw[<-] (-3,1.1) -- (-2.5,1.1);
\draw[<-] (-2.5,-0.9) -- (-2,-0.9);
\draw[<-] (-2,1.5) -- (-1.5,1.5);
\draw[<-] (-1.5,-1.2) -- (-1,-1.2);
\draw[<-] (-1,0.32) -- (-0.5,0.32);
\draw[<-] (-0.5,-0.32) -- (0,-0.32);
\draw[->] (0,-1.3) -- (0.5,-1.3);
\draw[->] (0,0.2) -- (0.5,0.2);
\draw[->] (1,-1.86) -- (1.5,-1.86);
\draw[->] (1.5,-1) -- (2,-1);
\draw[->] (2,0.7) -- (2.5,0.7);
\draw[<-] (3,-1.5) -- (3.5,-1.5);

% marcas contidos/steaflers 
\draw[->,line width=0.35mm] (-3.5,0.75) -- (-3,0.75);
\draw[->,line width=0.35mm] (-3,-1.8) -- (-2.5,-1.8);
\draw[<-,line width=0.35mm] (-2,-0.2) -- (-1.5,-0.2);
\draw[<-,line width=0.35mm] (-1,1.8) -- (-0.5,1.8);
\draw[->,line width=0.35mm] (-0.5,-0.7) -- (0,-0.7);
\draw[->,line width=0.35mm] (0,1.3) -- (0.5,1.3);
\draw[->,line width=0.35mm] (1,1.35) -- (1.5,1.35);
\draw[->,line width=0.35mm] (2.5,-0.25) -- (3,-0.25);
\draw[<-,line width=0.35mm] (3,0.72) -- (3.5,0.72);
\draw[->,line width=0.35mm] (-1.5,1) -- (-1,1);

% marcas para virar ignorantes/ignorants
\draw plot[mark=x] (-3.5,-1.1);
\draw plot[mark=x] (-3.5,1.15);
\draw plot[mark=x] (-2.5,-0.4);
\draw plot[mark=x] (-1.5,1.2);
\draw plot[mark=x] (0,1.55);
\draw plot[mark=x] (1,-1.7);
\draw plot[mark=x] (1,0.27);
\draw plot[mark=x] (2,1.7);
\draw plot[mark=x] (2.5,-1.22);
\draw plot[mark=x] (2.5,1.06);
\draw plot[mark=x] (3.5,1.52); 

% condições iniciais

\draw (-3.5,-2) node[below,font=\footnotesize] {$0$};
\draw (-3,-2) node[below,font=\footnotesize] {$2$};
\draw (-2.5,-2) node[below,font=\footnotesize] {$0$};
\draw (-2,-2) node[below,font=\footnotesize] {$1$};
\draw (-1.5,-2) node[below,font=\footnotesize] {$2$};
\draw (-1,-2) node[below,font=\footnotesize] {$0$};
\draw (-0.5,-2) node[below,font=\footnotesize] {$0$};
\draw (0,-2) node[below,font=\footnotesize] {$1$};
\draw (0.5,-2) node[below,font=\footnotesize] {$1$};
\draw (1,-2) node[below,font=\footnotesize] {$1$};
\draw (1.5,-2) node[below,font=\footnotesize] {$0$};
\draw (2,-2) node[below,font=\footnotesize] {$2$};
\draw (2.5,-2) node[below,font=\footnotesize] {$1$};
\draw (3,-2) node[below,font=\footnotesize] {$0$};
\draw (3.5,-2) node[below,font=\footnotesize] {$1$};

% evolução informantes

\draw[thick] (-2,-2) -- (-2,2);
\draw[thick] (-2.5,-0.9) -- (-2.5,-0.4);
\draw[thick] (-1,0.32) -- (-1,1.8);
\draw[thick] (-0.5,-0.32) -- (-0.5,2);
\draw[thick] (0.5,-2) -- (0.5,2);
\draw[thick] (0,-2) -- (0,1.55);
\draw[thick] (1,-2) -- (1,-1.7);
\draw[thick] (1.5,-1.86) -- (1.5,2);
\draw[thick] (2.5,-2) -- (2.5,-1.22);
\draw[thick] (2.5,0.7) -- (2.5,1.06);
\draw[thick] (3,-1.5) -- (3,2);
\draw[thick] (3.5,-2) -- (3.5,1.52);
\draw[thick] (-2.5,0) -- (-2.5,2);

% evolução adaptors

\draw[line width=0.5mm] (-3,-2) -- (-3,2);
\draw[line width=0.5mm] (-1.5,-2) -- (-1.5,1.2);
\draw[line width=0.5mm] (-1,1.8) -- (-1,2);
\draw[line width=0.5mm] (-1,1) -- (-1,2);
\draw[line width=0.5mm] (-2,-0.2) -- (-2,2);
\draw[line width=0.5mm] (2,-2) -- (2,1.7);

% evolução ignorants

\draw[dashed] (-3.5,-2) -- (-3.5,2);
\draw[dashed] (-2.5,-2) -- (-2.5,-0.9);
\draw[dashed] (-2.5,-0.4) -- (-2.5,2);
\draw[dashed] (-1.5,1.2) -- (-1.5,2);
\draw[dashed] (-1,-2) -- (-1,0.32);
\draw[dashed] (-0.5,-2) -- (-0.5,-0.32);
\draw[dashed] (0,1.55) -- (0,2);
\draw[dashed] (1,-1.7) -- (1,2);
\draw[dashed] (1.5,-2) -- (1.5,-1.86);
\draw[dashed] (2,1.7) -- (2,2);
\draw[dashed] (2.5,-1.22) -- (2.5,2);
\draw[dashed] (3,-2) -- (3,-1.5);
\draw[dashed] (3.5,1.52) -- (3.5,2);

\end{tikzpicture}
\end{center}
\caption{Realization of the graphical representation of the innovation diffusion model on $\mathbb{Z}$. The ignorant, aware and adopter states are drawn in dashed, thin and thick lines, respectively. {\color{black}The
thin arrows indicate the times at which an ignorant becomes aware, the thick arrows the times at
which an aware becomes adopter, and the $\times$ marks the times at which an agent becomes ignorant.}}
\end{figure}

%%%%%%%%%%%%%%%%%%%%%%%%%%%%%%%%%%%%%%%%%%%%%%%%%%%%
%%%%%%%% SS: Behavior of the Innovation Process
%%%%%%%%%%%%%%%%%%%%%%%%%%%%%%%%%%%%%%%%%%%%%%%%%%%%

\subsection{Behavior of the Innovation Process}

Consider the innovation process on $\mathbb{Z}^d$ with rates $\lambda$ and $\alpha$, for $d \geq 1$. We focus our attention on the influence of the parameters on the qualitative behavior of the process. First, we wish to obtain sufficient conditions under which the innovation awareness either becomes extinct or succeeds. In what follows, extinction has two different meanings depending on whether the initial configuration has finitely or infinitely many aware and adopters.

\begin{definition}
{\bf Extinction of the innovation awareness.} Whenever the initial configuration has finitely many aware and adopters, the innovation awareness is said to become extinct if there is almost surely a finite random time after which all sites in $\mathbb{Z}^d$ are ignorants (i.e., in state $0$). If the initial configuration has infinitely many aware or adopters, the innovation awareness is said to become extinct if for any fixed site there is almost surely a finite random time after which the site will stay in state $0$ forever. If the innovation awareness does not become extinct, we say that it is successful.
\end{definition}

Our first result states that the extinction or not of the innovation awareness, {\color{black}in the two senses defined above}, just depends on the value of $\lambda$. Let $\lambda_c(d)$ be the critical value of the basic $d$-dimensional contact process.  

\begin{theorem}\label{awareness}
{\color{black}For $d\geq 1$,} the innovation awareness becomes extinct if, and only if, $\lambda \leq \lambda_c(d)$.
\end{theorem}

{\color{black}
The previous theorem claims that the innovation process exhibits a phase transition in $\lambda$ regarding the extinction or success of the innovation awareness, and that the critical parameter coincides with the critical value for survival for the contact process, for any initial configuration. In the sequel, we focus our attention to the spread of the innovation adoption.}

\begin{definition}
{\bf Extinction of the innovation adoption.} If the initial configuration has finitely many adopters the innovation adoption is said to become extinct if there is almost surely a finite random time after which all sites in $\mathbb{Z}^d$ are ignorants or aware (i.e., in state $0$ or $1$). If the initial configuration has infinitely many adopters the innovation adoption is said to become extinct if for any fixed site there is almost surely a finite random time after which the site will stay in state $0$ or $1$ forever. As in the previous definition, if the innovation adoption does not become extinct, we say that it is successful.
\end{definition}

By Theorem \ref{awareness} we have that the innovation awareness is successful provided the rate $\lambda$ is greater than {\color{black}the} critical value of the basic $d$-dimensional contact process. However, this is not enough to guarantee the {\color{black}success of the innovation adoption. In the next result we show that a different phase transition appears when we analyze the extinction of the innovation adoption provided that $\lambda > \lambda_c(d)$.}

{\color{black}
\begin{theorem}\label{adoption}
For $d\geq 1$ and $\lambda > \lambda_c(d)$, there exists a critical value of the parameter $\alpha$, denoted $\alpha_c$, $\alpha_c\in (0,\infty)$, such that for any initial configuration with finitely many adopters
\begin{enumerate}
\item[(i)]  the innovation adoption becomes extinct, if $\alpha < \alpha_c$;
\item[(ii)] the innovation adoption is successful, if $\alpha > \alpha_c$.
\end{enumerate}
The same result holds for any initial configuration with infinitely many adopters for a possibly different critical value $\alpha_c$.
\end{theorem}}

\begin{remark}
Let us discuss how large is the domain where the statements of Theorems \ref{awareness} and \ref{adoption} hold. Observe that these results are stated invoking the critical parameter $\lambda_c(d)$. Although, up to now, it was not possible to obtain rigorous results regarding the exact numerical value of this parameter, some upper and lower {\color{black}bounds} have been obtained. A comparison between the contact process  and a branching random walk provides the lower bound $\lambda_c(d) \geq 1/(2d)$. On the other hand, Liggett \cite{improvedupperbound} proved that the critical parameter of the basic contact process, e.g. the contact process on $\mathbb{Z}, \lambda_c(1) \leq 1.9442$.
\end{remark}

%%%%%%%%%%%%%%%%%%%%%%%%%%%%%%%%%%%%%%%%%%%%%%%%%%%%
%%%%%%%% S: PROOFS
%%%%%%%%%%%%%%%%%%%%%%%%%%%%%%%%%%%%%%%%%%%%%%%%%%%%

\section{Proofs}
\label{sec:proofs}

This section is entirely devoted to the proof of our results. The approach pursued here is constructive and the main ideas are supported by the graphical representation of the innovation process. The idea behind  the proof of Theorem \ref{awareness} is a comparison between the innovation process and the basic contact process. The main argument used to prove Theorem \ref{adoption} is to compare the original process {\color{black}properly rescaled in space and time} with suitable oriented percolation models. The comparison of spatial stochastic processes with oriented percolation models is a powerful technique in the analysis of stochastic growth models --- {\color{black}see for instance \cite[Section 4]{Durret} where this technique is explained in the general context. See also \cite{raey,Schinazi2004}, and references therein, for some particular applications.}

\begin{proof}
{\bf Proof of Theorem \ref{awareness}.} We consider a coupling between the innovation process  on $\mathbb{Z}^d$, with rates $\lambda$ and $\alpha$, and the basic $d$-dimensional contact process $(\xi_t)_{t\geq 0}$ with infection rate $\lambda$ and death rate $1$. That is, if the process $(\xi_t)_{t \geq 0}$ is in configuration $\xi \in \{0,1\}^{\mathbb{Z}^d}$, the state of site $x$ changes according to the following transition rates

\begin{equation}\label{ratescontact}
\begin{array}{rclc}
0 & \rightarrow & 1, & \hspace{.5cm} \lambda \, n_1(x,\xi),\\
1 & \rightarrow & 0, & \hspace{.5cm}1.\\
\end{array}
\end{equation}

Indeed, the contact process is constructed, from the innovation process, as a new process which does not differentiate the states $1$ and $2$. At time $0$ we set $\xi_0(x)=1$ if $\eta_0(x)\neq 0$ and $\xi_0(x)=0$ if $\eta_0(x) =0$. In other words, we use as initial configuration for the contact process the initial configuration of the innovation process replacing all the $2$'s by $1$'s. Then we obtain a version of the contact process  with the rates given by {\color{black} \eqref{ratescontact}} using the same Poisson processes $N^{x,y}_{\lambda}$ and $D^x$ in the graphical construction of the innovation process and ignoring the marks of the Poisson process $N^{x,y}_{\alpha}$ (see Fig. 2). Thus defined, it is not difficult to see that the event of survival for the contact process is equivalent to the {\color{black}success} of innovation awareness in our model. Therefore, the theorem is a consequence of this comparison, the definition of the critical value $\lambda_c(d)$, {\color{black}  and the fact that the
critical contact process dies out (see \cite{bezuidenhout1990})}.

\end{proof}

\begin{figure}[h]\label{comparCP_SSIP}
\begin{center}

\begin{center}
\begin{tikzpicture}[scale=0.8]
\draw (-4.5,-2) -- (4.5,-2);
\draw (-4,-2) -- (-4,2);
\draw (4,-2) -- (4,2);

% linha de tempo
\draw[->] (-5,-2) -- (-5,2.5)
node [left,text width=2cm, midway]
{Innovation\\ Process};

% marcas informantes/spreaders
\draw[->,line width=0.35mm] (-3,0) -- (-2.5,0);
\draw[<-,line width=0.35mm] (-3,1.1) -- (-2.5,1.1);
\draw[<-,line width=0.35mm] (-2.5,-0.9) -- (-2,-0.9);
\draw[<-,line width=0.35mm] (-2,1.5) -- (-1.5,1.5);
\draw[->,line width=0.35mm] (-1.5,1) -- (-1,1);
\draw[<-,line width=0.35mm] (-1.5,-1.2) -- (-1,-1.2);
\draw[<-,line width=0.35mm] (-1,0.32) -- (-0.5,0.32);
\draw[<-,line width=0.35mm] (-0.5,-0.32) -- (0,-0.32);
\draw[->,line width=0.35mm] (0,-1.3) -- (0.5,-1.3);
\draw[->,line width=0.35mm] (0,0.2) -- (0.5,0.2);
\draw[->,line width=0.35mm] (1,-1.86) -- (1.5,-1.86);
\draw[->,line width=0.35mm] (1.5,-1) -- (2,-1);
\draw[->,line width=0.35mm] (2,0.7) -- (2.5,0.7);
\draw[<-,line width=0.35mm] (3,-1.5) -- (3.5,-1.5);

% marcas contidos/steaflers 
\draw[->,line width=0.5mm] (-3.5,0.75) -- (-3,0.75);
\draw[->,line width=0.5mm] (-3,-1.8) -- (-2.5,-1.8);
\draw[<-,line width=0.5mm] (-2,-0.2) -- (-1.5,-0.2);
\draw[<-,line width=0.5mm] (-1,1.8) -- (-0.5,1.8);
\draw[->,line width=0.5mm] (-0.5,-0.7) -- (0,-0.7);
\draw[->,line width=0.5mm] (0,1.3) -- (0.5,1.3);
\draw[->,line width=0.5mm] (1,1.35) -- (1.5,1.35);
\draw[->,line width=0.5mm] (2.5,-0.25) -- (3,-0.25);
\draw[<-,line width=0.5mm] (3,0.72) -- (3.5,0.72);

% marcas para virar ignorantes/ignorants
\draw plot[mark=x] (-3.5,-1.1);
\draw plot[mark=x] (-3.5,1.15);
\draw plot[mark=x] (-2.5,-0.4);
\draw plot[mark=x] (-1.5,1.2);
\draw plot[mark=x] (0,1.55);
\draw plot[mark=x] (1,-1.7);
\draw plot[mark=x] (1,0.27);
\draw plot[mark=x] (2,1.7);
\draw plot[mark=x] (2.5,-1.22);
\draw plot[mark=x] (2.5,1.06);
\draw plot[mark=x] (3.5,1.52); 

% condições iniciais

\draw (-3.5,-2) node[below,font=\footnotesize] {$0$};
\draw (-3,-2) node[below,font=\footnotesize] {$2$};
\draw (-2.5,-2) node[below,font=\footnotesize] {$0$};
\draw (-2,-2) node[below,font=\footnotesize] {$1$};
\draw (-1.5,-2) node[below,font=\footnotesize] {$2$};
\draw (-1,-2) node[below,font=\footnotesize] {$0$};
\draw (-0.5,-2) node[below,font=\footnotesize] {$0$};
\draw (0,-2) node[below,font=\footnotesize] {$1$};
\draw (0.5,-2) node[below,font=\footnotesize] {$1$};
\draw (1,-2) node[below,font=\footnotesize] {$1$};
\draw (1.5,-2) node[below,font=\footnotesize] {$0$};
\draw (2,-2) node[below,font=\footnotesize] {$2$};
\draw (2.5,-2) node[below,font=\footnotesize] {$1$};
\draw (3,-2) node[below,font=\footnotesize] {$0$};
\draw (3.5,-2) node[below,font=\footnotesize] {$1$};

% transmissão informantes

\draw[thick,line width=0.35mm] (-2,-2) -- (-2,2);
\draw[thick,line width=0.35mm] (-2.5,-0.9) -- (-2.5,-0.4);
\draw[thick,line width=0.35mm] (-1,0.32) -- (-1,1.8);
\draw[thick,line width=0.35mm] (-0.5,-0.32) -- (-0.5,2);
\draw[thick,line width=0.35mm] (0.5,-2) -- (0.5,1.3);
\draw[thick,line width=0.35mm] (0,-2) -- (0,1.55);
\draw[thick,line width=0.35mm] (1,-2) -- (1,-1.7);
\draw[thick,line width=0.35mm] (1.5,-1.86) -- (1.5,2);
\draw[thick,line width=0.35mm] (2.5,-2) -- (2.5,-1.22);
\draw[thick,line width=0.35mm] (3,-1.5) -- (3,0.72);
\draw[thick,line width=0.35mm] (3.5,-2) -- (3.5,1.52);

% evolução contidos

\draw[thick,line width=0.5mm] (-3,-2) -- (-3,2);
\draw[thick,line width=0.5mm] (-1.5,-2) -- (-1.5,1.2);
\draw[thick,line width=0.5mm] (-1,1.8) -- (-1,2);
\draw[thick,line width=0.5mm] (0.5,1.3) -- (0.5,2);
\draw[thick,line width=0.5mm] (2,-2) -- (2,1.7);
\draw[thick,line width=0.5mm] (3,0.72) -- (3,2);

% evolução ignorantes

\draw[dashed] (-3.5,-2) -- (-3.5,2);
\draw[dashed] (-2.5,-2) -- (-2.5,-0.9);
\draw[dashed] (-2.5,-0.4) -- (-2.5,2);
\draw[dashed] (-1.5,1.2) -- (-1.5,2);
\draw[dashed] (-1,-2) -- (-1,0.32);
\draw[dashed] (-0.5,-2) -- (-0.5,-0.32);
\draw[dashed] (0,1.55) -- (0,2);
\draw[dashed] (1,-1.7) -- (1,2);
\draw[dashed] (1.5,-2) -- (1.5,-1.86);
\draw[dashed] (2,1.7) -- (2,2);
\draw[dashed] (2.5,-1.22) -- (2.5,2);
\draw[dashed] (3,-2) -- (3,-1.5);
\draw[dashed] (3.5,1.52) -- (3.5,2);

\end{tikzpicture}
\end{center}

\begin{center}
\begin{tikzpicture}[scale=0.8]
\draw (-4.5,-2) -- (4.5,-2);
\draw (-4,-2) -- (-4,2);
\draw (4,-2) -- (4,2);

% linha de tempo
\draw[->] (-5,-2) -- (-5,2.5)
node [left,text width=2cm, midway]
{Contact\\ process};

% marcas informantes/spreaders
\draw[->,line width=0.35mm] (-3,0) -- (-2.5,0);
\draw[<-,line width=0.35mm] (-3,1.1) -- (-2.5,1.1);
\draw[<-,line width=0.35mm] (-2.5,-0.9) -- (-2,-0.9);
\draw[<-,line width=0.35mm] (-2,1.5) -- (-1.5,1.5);
\draw[->,line width=0.35mm] (-1.5,1) -- (-1,1);
\draw[<-,line width=0.35mm] (-1.5,-1.2) -- (-1,-1.2);
\draw[<-,line width=0.35mm] (-1,0.32) -- (-0.5,0.32);
\draw[<-,line width=0.35mm] (-0.5,-0.32) -- (0,-0.32);
\draw[->,line width=0.35mm] (0,-1.3) -- (0.5,-1.3);
\draw[->,line width=0.35mm] (0,0.2) -- (0.5,0.2);
\draw[->,line width=0.35mm] (1,-1.86) -- (1.5,-1.86);
\draw[->,line width=0.35mm] (1.5,-1) -- (2,-1);
\draw[->,line width=0.35mm] (2,0.7) -- (2.5,0.7);
\draw[<-,line width=0.35mm] (3,-1.5) -- (3.5,-1.5);

% marcas contidos/steaflers 
\draw[->,line width=0.5mm] (-3.5,0.75) -- (-3,0.75);
\draw[->,line width=0.5mm] (-3,-1.8) -- (-2.5,-1.8);
\draw[<-,line width=0.5mm] (-2,-0.2) -- (-1.5,-0.2);
\draw[<-,line width=0.5mm] (-1,1.8) -- (-0.5,1.8);
\draw[->,line width=0.5mm] (-0.5,-0.7) -- (0,-0.7);
\draw[->,line width=0.5mm] (0,1.3) -- (0.5,1.3);
\draw[->,line width=0.5mm] (1,1.35) -- (1.5,1.35);
\draw[->,line width=0.5mm] (2.5,-0.25) -- (3,-0.25);
\draw[<-,line width=0.5mm] (3,0.72) -- (3.5,0.72);

% marcas para virar ignorantes/ignorants
\draw plot[mark=x] (-3.5,-1.1);
\draw plot[mark=x] (-3.5,1.15);
\draw plot[mark=x] (-2.5,-0.4);
\draw plot[mark=x] (-1.5,1.2);
\draw plot[mark=x] (0,1.55);
\draw plot[mark=x] (1,-1.7);
\draw plot[mark=x] (1,0.27);
\draw plot[mark=x] (2,1.7);
\draw plot[mark=x] (2.5,-1.22);
\draw plot[mark=x] (2.5,1.06);
\draw plot[mark=x] (3.5,1.52); 

% condições iniciais

\draw (-3.5,-2) node[below,font=\footnotesize] {$0$};
\draw (-3,-2) node[below,font=\footnotesize] {$1$};
\draw (-2.5,-2) node[below,font=\footnotesize] {$0$};
\draw (-2,-2) node[below,font=\footnotesize] {$1$};
\draw (-1.5,-2) node[below,font=\footnotesize] {$1$};
\draw (-1,-2) node[below,font=\footnotesize] {$0$};
\draw (-0.5,-2) node[below,font=\footnotesize] {$0$};
\draw (0,-2) node[below,font=\footnotesize] {$1$};
\draw (0.5,-2) node[below,font=\footnotesize] {$1$};
\draw (1,-2) node[below,font=\footnotesize] {$1$};
\draw (1.5,-2) node[below,font=\footnotesize] {$0$};
\draw (2,-2) node[below,font=\footnotesize] {$1$};
\draw (2.5,-2) node[below,font=\footnotesize] {$1$};
\draw (3,-2) node[below,font=\footnotesize] {$0$};
\draw (3.5,-2) node[below,font=\footnotesize] {$1$};

% transmissão informantes

\draw[thick,line width=0.35mm] (-2,-2) -- (-2,2);
\draw[thick,line width=0.35mm] (-2.5,-0.9) -- (-2.5,-0.4);
\draw[thick,line width=0.35mm] (-1,0.32) -- (-1,1.8);
\draw[thick,line width=0.35mm] (-0.5,-0.32) -- (-0.5,2);
\draw[thick,line width=0.35mm] (0.5,-2) -- (0.5,1.3);
\draw[thick,line width=0.35mm] (0,-2) -- (0,1.55);
\draw[thick,line width=0.35mm] (1,-2) -- (1,-1.7);
\draw[thick,line width=0.35mm] (1.5,-1.86) -- (1.5,2);
\draw[thick,line width=0.35mm] (2.5,-2) -- (2.5,-1.22);
\draw[thick,line width=0.35mm] (3,-1.5) -- (3,0.72);
\draw[thick,line width=0.35mm] (3.5,-2) -- (3.5,1.52);

% evolução contidos

\draw[thick,line width=0.35mm] (-3,-2) -- (-3,2);
\draw[thick,line width=0.35mm] (-1.5,-2) -- (-1.5,1.2);
\draw[thick,line width=0.35mm] (-1,1.8) -- (-1,2);
\draw[thick,line width=0.35mm] (0.5,1.3) -- (0.5,2);
\draw[thick,line width=0.35mm] (2,-2) -- (2,1.7);
\draw[thick,line width=0.35mm] (3,0.72) -- (3,2);

% evolução ignorantes

\draw[dashed] (-3.5,-2) -- (-3.5,2);
\draw[dashed] (-2.5,-2) -- (-2.5,-0.9);
\draw[dashed] (-2.5,-0.4) -- (-2.5,2);
\draw[dashed] (-1.5,1.2) -- (-1.5,2);
\draw[dashed] (-1,-2) -- (-1,0.32);
\draw[dashed] (-0.5,-2) -- (-0.5,-0.32);
\draw[dashed] (0,1.55) -- (0,2);
\draw[dashed] (1,-1.7) -- (1,2);
\draw[dashed] (1.5,-2) -- (1.5,-1.86);
\draw[dashed] (2,1.7) -- (2,2);
\draw[dashed] (2.5,-1.22) -- (2.5,2);
\draw[dashed] (3,-2) -- (3,-1.5);
\draw[dashed] (3.5,1.52) -- (3.5,2);

\end{tikzpicture}
\end{center}

\caption{Coupling between the innovation process and the basic contact process. For any $x\in \mathbb{Z}^d$ and any $t\geq 0$, if $\eta_t(x)=1$ or $\eta_t(x)=2$ then $\xi_t(x)=1$. {\color{black} The solid lines in the bottom picture are the sites in state $1$ for the contact process.}}
\end{center}
\end{figure}
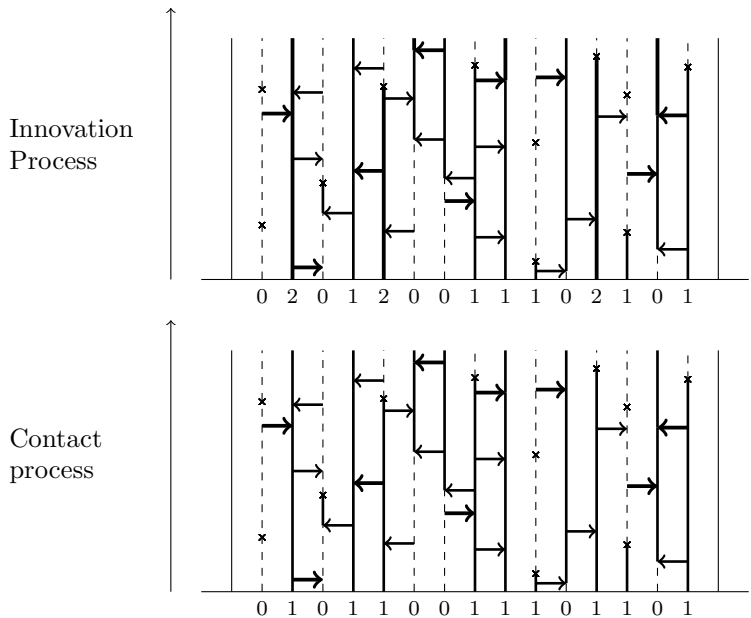

\begin{proof}
{\bf Proof of Theorem \ref{adoption}.}

{\color{black}   For $d\geq 1$, let $\mathcal{A}_0:=\{x\in \mathbb{Z}^d: \eta_0(x)=2\}$ and
$$\theta(\lambda,\alpha,d):=\mathbb{P}\left[\bigcup_{x\in\mathbb{Z}^d}\{\eta_t(x)=2\},\text{ for all }t>0 \Bigg| \mathcal{A}_0=\{{\bf 0}\}\right].$$
We observe that the innovation adoption becomes extinct if, and only if, $\theta(\lambda,\alpha,d)=0$.  In order to prove Theorem \ref{adoption}, we first observe that the function $\theta(\alpha):=\theta(\lambda,\alpha,d)$ is nondecreasing in $\alpha\in (0,\infty)$. We can verify this claim by using the graphical construction of the innovation model. Indeed, we can construct two innovation models with rates $\lambda, \alpha_1$, and $\lambda, \alpha_2$, respectively with $\alpha_1 < \alpha_2$ and starting from the same initial configuration. Then, both models can be coupled in such a way that whenever a site is in state $2$ for the process associated to the rates $\lambda, \alpha_1$, then the same is true for the process with rates $\lambda, \alpha_2$. This in turn implies that the success of innovation adoption for the first model implies the success of the innovation adoption for the second model too. For the sake of simplicity, we assume that $\mathcal{A}_0=\{{\bf 0}\}$ in the definition of $\theta(\alpha,\lambda,d)$. However, the same result about monotonicity is true for any initial configuration with at least one site in state $2$. The reader may see in \cite[Section 6.2]{Grimmett} an application of the coupling method and the graphical construction for the study of the monotonicity of the contact process. The previous remark implies that the critical value

\begin{equation}\label{criticalalpha}
\alpha_c :=\alpha_c(\lambda,d)=\sup\left\{\alpha \in[0,\infty):\text{the innovation adoption becomes extinct}\right\}
\end{equation}

is well defined. {\color{black}Thus}, the crucial point in the proof of Theorem \ref{adoption} is showing that $\alpha_c \in (0,\infty)$, which is a consequence of the following propositions, stated for $d\geq 1$ and $\lambda > \lambda_c(d)$, and for any initial configuration starting with at least one site in state $2$.

\begin{prop}\label{prop1}
For $\alpha$ sufficiently close to $0$, the innovation adoption becomes extinct.
\end{prop}

\begin{prop}\label{prop2}
For $\alpha$ sufficiently large, the innovation adoption is successful.
\end{prop}

We prove these results below.
}

\subsubsection*{Proof of Prop. \ref{prop1}: {\bf Extinction of the innovation adoption.}} The main idea in order to prove that the innovation adoption becomes extinct is to compare the process $(\eta_t)_{t\geq 0}$ with a suitable oriented percolation model defined on $\mathcal{L}:=\mathbb{Z}^d \times \mathbb{Z}^+$. The first step in this direction is to consider the nested space-time regions (see Fig. \ref{caixapb})
\begin{eqnarray}
\Lambda_1 = [-2L,2L]^d \times [0,2T] \ \ \ \ \ \ \ \Lambda_2 = [-L,L]^d \times [T,2T].
\end{eqnarray}
Let $\Delta$ be the boundary of $\Lambda_1$, defined as
\[
\Delta = \{(x,m) \in \Lambda_1 : |x_i|=2L \ \mbox{for some} \ i, i=1, \ldots d \ \mbox{or} \ m=0\}.
\]

\begin{figure}[!ht]
	\centering
	\includegraphics[width=0.85\linewidth]{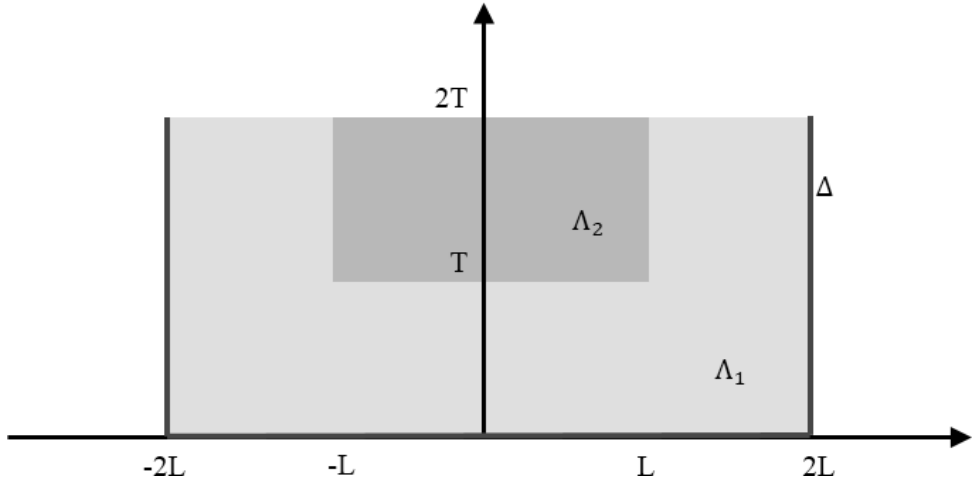}
	\caption{Schema of the boxes for $d = 1$, where $\mathcal{L} = \mathbb{Z} \times \mathbb{Z}_{+}$, $\Lambda_{1} = [-2L,2L]\times [0,2T]$ and $\Lambda_{2} = [-L,L]\times [T,2T]$.}
	\label{caixapb}
\end{figure}

A site percolation model is defined on the lattice $\mathcal{L}$ by declaring a site $(x,m) \in \mathcal{L}$ as open if, and only if, for the process $(\eta_t)_{t\geq 0}$ restricted to $\Lambda_1 + (x,m)$ the box $\Lambda_2 + (x,m)$ has {\color{black}no} sites in state $2$ regardless of the states of sites in the boundary $\Delta + (x,m)$. As usual, sites which are not open are called closed. In this sense, if we consider for each $(x,m)\in \mathcal{L}$ the random variable 
$$X_{(x,m)}=\left\{
\begin{array}{ll}
1, & \textrm{ if }(x,m)\textrm{ is open,}\\
0,& \textrm{ otherwise}
\end{array}\right.
$$ 
then, the collection of random variables $\mathcal{X}:=\{X_{(x,m)}:(x,m)\in \mathcal{L}\}$ induces the desired percolation model on $\mathcal{L}$. To see this, we first make $\mathcal{L}$ into a {\color{black} directed} graph. Denote $\Lambda(x,m) := \Lambda_{1} + (x,m)$ and draw an oriented edge from $(x,m)$ to $(y,n)$ if, and only if, $m \leq n$ and $\Lambda(x,m) \cap \Lambda(y,n) \neq \emptyset$. The open sites in the resulting directed graph define a percolation model as follows. We say that $(y,n)$ can be reached from $(x,m)$  and write $(x,m) \rightarrow (y,n)$ if there is a sequence of sites $x_{0} = x,...,x_{k} = y$ and time instants $n_{0} = m,...,n_{k} =n$ such that: first, there is an oriented edge from $(x_{i},n_{i})$ to $(x_{i+1},n_{i+1})$ for $0 \leq i < k$; second, $X_{(x_{i},n_{i}) }= 1$ for $0 \leq i \leq k$. The resulting model is a $K$-dependent percolation model on $\mathcal{L}$. Indeed, there exists a constant $K$ depending only on the dimension $d$ such that if the distance between sites $(x,m)$ and $(y,n)$ is larger than $K$ then the associated random variables $X_{(x,m)}$ and $X_{(y,n)}$ are independent. 

Now, the main idea is to show that for any $\varepsilon > 0$, and for the induced percolation model, 
{\color{black}if $\alpha$ is sufficiently close to $0$, then}
	\begin{equation}\label{open}
		\mathbb{P}[(0,0) \text{ is open}]\geq 1 - 2 \varepsilon.
	\end{equation}
	
Note that by translation-invariance of the process, \eqref{open} implies that 
$$\mathbb{P}[(x,n) \text{ is open}]\geq 1 - 2 \varepsilon$$ 
for any $(x,n)\in \mathcal{L}$. Once we have \eqref{open}, the rest of the proof is somewhat standard and we refer the reader to Van Den Berg et al. \cite{vandenberg1998} for more details. The crucial point is that if, for any $(x,m)$, there is an adopter in $\Lambda_2 + (x,m)$, for the innovation process then,    {\color{black} site $(x,m)$ can be reached from a path of closed sites in the associated percolation model}. By taking $\varepsilon$ small enough it is possible to make the probability of a path of closed sites decay exponentially fast with its length. This in turn implies that for any fixed site  in the innovation process there is a finite random time after which the site is in state $0$ or $1$. Hence, we have the extinction of the innovation adoption.

In order to prove \eqref{open}, suppose that $(0,0)$ is closed. That is, there is at least one site in state $2$ inside $\Lambda_{2}$. Define the event $A$ as being the event ``there are {\color{black}no} arrows of the process  $N_{\alpha}^{x,y}$ inside $\Lambda_{1}$ ,'' and note that
$$\mathbb{P}[(0,0) \text{ is open}] {\color{black}\geq} \mathbb{P}[(0,0) \text{ is open} | A]\cdot P(A).$$ 

Now, let $\varepsilon >0$. Observe that, conditioned on $A$, if there is one site in state $2$ in $\Lambda_{2}$ then, such site must be on that state from the lower part of the boundary $\Delta$, restricted to $[-L,L]^{d}$. Then, we have that
$$\mathbb{P}[(0,0) \text{ is closed} | A] \leq \mathbb{P}\left[\bigcup_{x \in [-L,L]^{d}} \{D^{x}(0,T) = 0\}\right] \leq (2L + 1)^{d}e^{-T}.$$
By taking $L = T$ large enough we get
$$
	\mathbb{P}[(0,0) \text{ is open} | A]  \geq 1 - (2L + 1)^{d}e^{-L} \geq 1 - \varepsilon.
$$
In addition, since $\Lambda_{1}$ is a finite space-time box we can pick $\alpha$ small enough such that the event $A$ occurs with probability at least $1 - \varepsilon$. Its proof may follows from the first moment method. Therefore, we conclude that for $\alpha$ small enough
$$
	\mathbb{P}[(0,0)  \text{ is open}] \geq  1 - 2 \varepsilon.
$$
This completes the proof of {\color{black}Proposition \ref{prop1}.}

\subsubsection*{Proof of Prop. \ref{prop2}: {\bf Survival of the innovation adoption.}}
As in the previous proof, we will compare the innovation process with a suitable oriented percolation model. Consider 
$$
	\mathcal{L}_{0} = \{(x,n) \in \mathbb{Z}^{2}: x + n \hspace{0.1cm} \text{is even;} \hspace{0.1cm} n \geq 0\}
$$
and
$$
	B = (-4L,4L)^{d} \times [0,T] \hspace{2.0cm} B_{x,n} = (2xL,nT) + B
$$

$$
	I = [-L,L]^{d} \hspace{2.0cm} I_{x} = 2xL + I
$$

\noindent where $L$ and $T$ are values to be defined later (see Fig. \ref{fig:caixa2}). Let $k = [\sqrt{L}]$, where $[a]$ denotes the integer part of $a$. We declare a site $(x,n) \in \mathcal{L}_{0}$ as open if, and only if, at time $nT$ there are no sites in state $1$ in $I_{x}$ and there are at least $k$ sites in state $2$ in $I_{x}$, and if at time $(n+1)T$ there are no sites in state $1$ in $I_{x-1}$ and $I_{x+1}$ and there are at least $k$ sites in state $2$ in each interval.  As usual, sites of $\mathcal{L}_{0}$ which are not open are called closed.

Now, if we consider for each $(x,n)\in \mathcal{L}_0$ the random variable 
$$Y_{(x,n)}=\left\{
\begin{array}{ll}
1, & \text{ if }(x,n)\text{ is open,}\\
0,& \text{ otherwise}
\end{array}\right.
$$ 
then, the collection of random variables $\mathcal{Y}:=\{Y_{(x,n)}:(x,n)\in \mathcal{L}_0\}$ induces a $K$-dependent percolation model on $\mathcal{L}_0$. Furthermore, {\color{black}an infinite cluster of open sites} in the induced percolation model implies the existence of adopters at all times in the innovation process. That is, the survival of the innovation adoption. Therefore, in order to prove survival of the innovation process we will prove that, for $\varepsilon > 0$, 
{\color{black}if we consider $\alpha$ sufficiently large, then}
\begin{equation}\label{open2}
	\mathbb{P}[(x,n) \text{ is open}] \geq 1 - 3 \varepsilon.
\end{equation}
Thus, since we choose $\varepsilon$ to be small enough, the existence of percolation of open sites is a consequence of well known results on $K$-dependent oriented percolation (see \cite[Section 4]{Durret} for instance). 
\begin{figure}[!ht]
	\centering
	\includegraphics[width=1.0\linewidth]{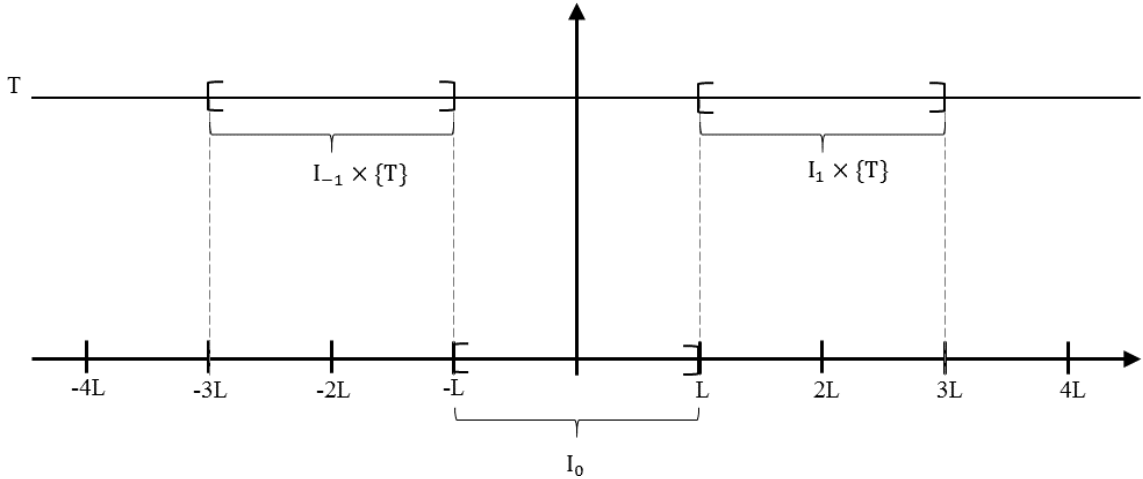}
	\caption{Schema for the site $(0,0)$ in the case $d = 1$.}
	\label{fig:caixa2}
\end{figure}

In order to prove \eqref{open2}, by translation invariance of the process, it suffices to prove it for site $(0,0)$. To simplify the notation we suppose in the sequel that $ d = 1$. However, our arguments hold for any $d \geq 1$. Let $\varepsilon >0$. Define the event $A_{x,y}$ as being the event ``every time that there is an arrival of the process $N_\lambda^{xy}$, before anything else happens {\color{black}in the box $B$}, there is an arrival of the process $N_\alpha^{xy}$ on $[0,T]$". Let

$$A =  \bigcap \limits_{{\stackrel{x, y \in (-4L,4L)}{||x - y|| = 1}}} A_{x,y},$$
and note that $\mathbb{P}[(x,n) \text{ is open}] {\color{black}\geq} \mathbb{P}[(x,n) \text{ is open}| A]\cdot P[A].$ It is not difficult to see that if we pick $\alpha$ large enough, then $\mathbb{P}(A) \geq 1 - \varepsilon$. Now, we will focus our attention on finding a lower bound for the probability of $(x,n)$ be open, conditioned to the occurrence of $A$.

We point out that by assuming that at time $0$ there are at least $k$ sites in state $2$ (adopters) and no sites in state $1$ (aware) in $I$, the adopters in $I$, conditioned on the event $A$, behave like a contact process in the following sense. Adopters survive at least as well as a contact process restricted to $I$. This is so because in our process adopters could appear from outside $I$ into such interval while this is not allowed in the contact process restricted to $I$. Now, consider the super-critical contact process restricted to the finite volume $\{1,2,...,k\}$ and let $\tau^{k}$ be the (random) time it takes for the restricted process to become extinct. It is a well known fact that 
$$
\mathbb{P}(\tau^{k} \leq e^{k}) < e^{-k} \hspace{0.2cm} \text{for large} \hspace{0.1cm} k,
$$

(see  Proposition 2.1 in \cite{Mountford(1993)}). In other words, this means that the supercritical contact process restricted to the finite volume $\{1,2,...,k\}$ survives at least $e^{k}$ with probability at least $1 - e^{-k}$. It is not difficult to show that this in turn implies that there exists a time $T_1$ such that the number of adopters in $I$ is, at that time, at least $M:=[\sqrt{k}]$ with probability at least
$$
(1 - e^{-M})^{M} \geq 1 - \varepsilon,
$$
for $T_1$ large enough. Indeed, $T_1$ may be taken as an increasing function of $L$. We refer the reader, for instance, to \cite{Schinazi2004} for a detailed proof of this claim. Thus at time $T_1$ we have at least $M$ adopters in $I$ and no aware in $[-4L + 1,4L - 1]$ with probability at least $1-\varepsilon$. Since $\lambda > \lambda_{c}(d)$, we can use well-known results of Bezuidenhout and Grimmett \cite{bezuidenhout1990} for the supercritical contact process. In particular, given that a supercritical contact process does not die out, the Shape Theorem (see Liggett \cite{opac-b1095541}, p 128) ensures that the adopters spread linearly. Hence, there is a constant $a > 0$ such that by time $aT_1$ the adopters have reached the sites $-4L$ and $4L$. Hence, by time $T_1 + aT_1$ there are at least $k$ adopters in $I_{1} = [L,3L]$ and in $I_{-1} = [-3L,-L]$ with probability at least $1 - 2 \varepsilon$ (one $\varepsilon$ takes care of the survival probability of $M$ adopters and the other one of the Shape Theorem). Set $T = T_1 + aT_1$ where $T_1$ is large enough, provided $L$ is large enough. Thus, we have shown that with probability at least $1 - 2 \varepsilon$ we will have at least $k$ adopters in $[-3L,-L]$ and in $[L,3L]$ at a certain time $T$.

{\color{black}Therefore}, if we pick $\alpha$ {\color{black}sufficiently large}, we can guarantee that

$$
\mathbb{P}[(0,0) \hspace{0.1cm} \textrm{is open} ] \geq \mathbb{P}[(0,0) \hspace{0.1cm} \textrm{is open}| A] \mathbb{P}(A) \geq (1-  2\varepsilon)(1 - \varepsilon) \geq 1 -3      \varepsilon.
$$
\noindent
This gives \eqref{open2} which finishes the proof.
\end{proof}

%%%%%%%%%%%%%%%%%%%%%%%%%%%%%%%%%%%%%%%%%%%%%%%%%%%%
%%%%%%%% S: CONCLUDING REMARKS
%%%%%%%%%%%%%%%%%%%%%%%%%%%%%%%%%%%%%%%%%%%%%%%%%%%%

\section{Concluding Remarks}\label{concluding}

In this paper we propose a simple mathematical model to capture the essence of the phenomenon of innovation diffusion on a structured population. In other words, we consider a model with a relationship between stochastic individual behavior and aggregate behavior. In this way, we complement recent studies regarding this issue. See for instance \cite{Goldenberg2004,Goldenberg2001,Goldenberg2002,Shun2002}. Our model is an interacting particle system and our results are obtained through the application of well known techniques of coupling between the innovation stochastic process, the contact process and suitable oriented percolation models. Such approach is an alternative to the techniques, growing in quantity and influence, used by various researchers in diffusion theory.  Moreover, our model is simple and the parameters have intuitive interpretations. These properties could lead to modifications, generalizations, and applications in more realistic scenarios. For example, consider the continuous-time Markov process $(\xi_t)_{t\geq 0}$ with {\color{black}state} space $\mathcal{S}=\{0,1,2\}^{\mathbb{Z}^d}$ and the same interpretation as before, e.g. $0$ representing an ignorant agent, $1$ an aware and $2$ an adopter. Then, if the system is in configuration $\xi \in \mathcal{S},$ assume that the state of site $x$ changes according to the following transition rates

\begin{equation*}\label{rates}
\begin{array}{rclc}
0 & \rightarrow & 1, & \hspace{.5cm} \lambda \, (n_1(x,\xi) + n_2(x,\xi)),\\ 
0 & \rightarrow & 2, & \hspace{.5cm}  \gamma,\\
1 & \rightarrow & 2, & \hspace{.5cm}\alpha \, n_2(x,\xi),\\
1 & \rightarrow & 0, & \hspace{.5cm}1,\\ 
2 & \rightarrow & 0, & \hspace{.5cm}1,\\ 
\end{array}
\end{equation*}

\noindent where $$n_i(x,\xi)= \sum_{||x-y||=1} 1\{\xi(y)=i\}$$ 
is the number of nearest neighbors of site $x$ in state $i$ for the configuration $\xi$, for $i=1,2.$ Thus defined, this model incorporates an innovation parameter as the Bass model does. {\color{black} This transition from state $0$ to state $2$ at rate $\gamma$ may represent the decision to adopt an innovation independently of the decisions of other individuals in a social system, and that could happen, for example, due to an external influence as mass media. Those individuals who exhibit this transition are called innovators in \cite{Bass1969}.}   We point out that by choosing $\gamma$ small enough, all our results also hold for this model. One can observe this by taking the sizes of the boxes, used in the proofs, such that with high probability there are no marks of the Poisson process with intensity $\gamma$ inside such boxes. The assumption that $\gamma$ has to be small is natural according to previous diffusion models or observed applications. See, for instance, the assumptions about the innovation parameter in \cite{Bass1969}.  

%%%%%%%%%%%%%%%%%%%%%%%%%%%%%%%%%%%%%%%%%%%%%%%%%%%%
%%%%%%%% S: ACKNOWLEDGEMENTS
%%%%%%%%%%%%%%%%%%%%%%%%%%%%%%%%%%%%%%%%%%%%%%%%%%%%

 \section*{Acknowledgements}
 
Special thanks are given to the referee, whose careful reading of the manuscript and valuable comments contributed to improve this paper. K.B.E.O. thanks FAPESP (Scholarship 12/22185-0), {\color{black}C.F.C. and P.M.R. thank CNPq (Grant 479313/2012-1), and P.M.R. also thanks FAPESP (Grants 2013/03898-8, 2015/03868-7), for financial support. The last author was visiting the Laboratoire de Probabilités et Modèles Aléatoires, Université Paris-Diderot when part of this work was carried out and he is grateful for their hospitality and support.}

%%%%%%%%%%%%%%%%%%%%%%%%%%%%%%%%%%%%%%%%%%%%%%%%%%%%
%%%%%%%% REFERENCES
%%%%%%%%%%%%%%%%%%%%%%%%%%%%%%%%%%%%%%%%%%%%%%%%%%%%

\end{document}